\documentclass[12pt]{article}
\usepackage{a4wide}

\newcommand{\matteop}[2]{\def#1{\mathop{\rm #2}\nolimits}}
\newcommand{\sstack}[2]{{\scriptstyle #1\atop\scriptstyle #2}}
\newcommand{\smallbox}[1]{\hbox to 0pt{\hss$\scriptstyle#1$\hss}}
\matteop{\Res}{Res}
\newcommand{\B}[1]{{\bf #1}}
\newcommand{\C}{\B{C}}
\newcommand{\Z}{\B{Z}}
\newcommand{\N}{\B{N}}
\renewcommand{\P}{\B{P}}
\newcommand{\pil}{\rightarrow}
\newcommand{\lpil}{\longrightarrow}

\newcommand{\rmap}[1]{\stackrel{#1}\lpil}

\matteop{\Pf}{Pf}
\newcommand{\SO}{{\cal O}}
\newcommand{\ds}{\displaystyle}
\matteop{\rank}{rank}

\def\Qterm#1#2{{\textstyle\ifnum#2=1\frac{#1}{1-#1}\else%
\frac{#2^3#1^#2}{1-#1^#2}\fi}}
\newcommand{\qterm}[1]{\Qterm{q}{#1}}
\newcommand{\qcterm}[1]{\ifnum#1=1\Qterm{\frac{q}{c_2}}{#1}\else
 \Qterm{(\frac{q}{c_2})}{#1}\fi}
\newcommand{\qtterm}[1]{\Qterm{\tilde q}{#1}}

\newtheorem{theorem}{Theorem}
\newtheorem{prop}[theorem]{Proposition}

\newtheorem{conj}[theorem]{Conjecture}
\newtheorem{defin}[theorem]{Definition}

\title{The Pfaffian Calabi--Yau, its Mirror, and their link to the
Grassmannian G(2,7)}

\author{Einar Andreas R\o{}dland\thanks{Dept.\ of Mathematics,
    University of Oslo, Box 1053 Blindern, 0316 Oslo, Norway; E-mail:
    {\tt einara@math.uio.no}; WWW: {\tt http://www.math.uio.no/$\sim$einara}.}}

\date{January 20, 1998\message{*** HUSK DATO ***}}

\begin{document}

\sloppy
\maketitle

\abstract{The rank 4 locus of a general skew-symmetric $7\times7$
matrix gives the pfaffian variety in $\P^{20}$ which is not defined as
a complete intersection. Intersecting this with
a general $\P^6$ gives a Calabi--Yau manifold. An orbifold
construction seems to give the 1-parameter mirror-family of this.
However, corresponding to two points in the 1-parameter family of
complex structures, both with maximally unipotent monodromy, are two
different mirror-maps: one corresponding to the general pfaffian
section, the other to a general intersection of $G(2,7)\subset\P^{20}$
with a $\P^{13}$. Apparently, the
pfaffian and $G(2,7)$ sections constitute different parts of the
A-model (K\"ahler structure related) moduli space, and, thus,
represent different parts of the same conformal field theory moduli
space.}

\section{The Pfaffian Variety}

Let $E$ be a rank 7 vector space. For $N\in E\wedge E$ non-zero, we
look at the locus of $\bigwedge^3N=0\in\bigwedge^6E$: the rank 4 locus of $N$
if viewed as a skew-symmetric matrix. This defines a degree 14
variety of codimension 3 in $\P(E\wedge E)\cong\P^{20}$. As $N$ is
skew-symmetric, this variety is defined by the pfaffians, ie.\ square
roots of the determinants, of the $6\times6$ diagonal minors of the
matrix.  Intersecting this with a general 6-plane in $\P(E\wedge
E)\cong\P^{20}$ will give a 3-dimensional Calabi--Yau. In coordinates
$x_i$ on $\P^6$, the matrix $N$ can be written $N_A=\sum_{i=0}^6
x_iA_i$ where the $A_i\in E\wedge E$ are skew-symmetric matrices
spanning the $\P^6$. Denote this variety $X_A\subset\P^{6}$. The
pfaffian variety in $\P^{20}$ is smooth away from the rank 2 locus
which has dimension 10. Hence, by Bertini's theorem, the variety $X_A$
is smooth for general $A$.

\begin{defin}
Let $N_A=\sum_{i=0}^6 x_iA_i$ where $A_i$ are $7\times7$ skew-symmetric
matrices. Let $X_A\subset\P^6$ denote the zero-locus of the pfaffians of the
$6\times6$ diagonal minors of $N_A$: ie., the rank 4 locus of the matrix.
\end{defin}

For $P=N^3\in\bigwedge^6E$, $\SO=\SO_{\P(E)}$, there are exact sequences
\begin{equation}
\begin{array}{c}
0\lpil (\bigwedge^7E^\vee)^2\otimes\SO(-7)
\rmap{P} \bigwedge^7E^\vee\otimes E^\vee\otimes\SO(-4)
\\
\quad
\rmap{N} \bigwedge^7E^\vee\otimes E\otimes\SO(-3)
\rmap{P} \SO
\lpil \SO_X
\lpil 0
\end{array}
\end{equation}
and
\begin{equation}
\begin{array}{c}
0\lpil (\bigwedge^7E^\vee)^2\otimes\bigwedge^2E^\vee\otimes\SO(-8)
\rmap{N\cdot}
 (\bigwedge^7E^\vee)^2\otimes({\rm Hom}(E,E)/{\rm Id}_E)\otimes\SO(-7)
\\
\quad
\rmap{\cdot N} (\bigwedge^7E^\vee)^2\otimes S^2E\otimes\SO(-6)
\rmap{P^{\otimes2}} {\cal J}_X^2
\lpil 0
\end{array}
\end{equation}
or more simply, for $P=[p_i]$ the pfaffians with proper choice of sign
and ordering,
\begin{equation}
0\lpil\SO_{\P^6}(-7)
\rmap{P^{\rm T}} 7\SO_{\P^6}(-4)
\rmap{N} 7\SO_{\P^6}(-3)
\rmap{P} \SO_{\P^6}
\lpil \SO_X
\lpil 0
\end{equation}
and
\begin{equation}
0 \lpil 21\SO_{\P^6}(-8) \lpil 48\SO_{\P^6}(-7)
\lpil 28\SO_{\P^6}(-6) \rmap{P^{\otimes2}} {\cal J}_X^2 \lpil 0.
\end{equation}
These sequences together with $0\pil{\cal J}_X^2\pil{\cal
  J}_X\pil{\cal N}_X^\vee\pil0$, $0\pil{\cal
  N}_X^\vee\pil\Omega_{\P^6}|_X\pil\Omega_X\pil 0$, and
$0\pil\Omega_{\P^6}|_X\pil7\SO_X(-1)\pil\SO_X\pil0$ give the
cohomology of the general, smooth manifold:

\begin{prop}
  The general variety $X_A\subset\P^6$ is smooth with
  $h^{1,0}=h^{2,0}=0$, $h^{3,0}=1$, $h^{1,1}=h^{2,2}=1$,
  $h^{1,2}=h^{2,1}=50$, $\chi=-98$, and $\omega_X\cong\mathop{\it{\cal
      E}\!xt}\nolimits^3(\SO_{X_A},\omega_{\P^6}) \cong\SO_{X_A}$;
  hence, it is a Calabi--Yau manifold. When $X_A$ is singular, we have
  trivial dualizing sheaf, $\omega_{X_A}^\circ\cong\SO_{X_A}$.
\end{prop}

\section{The Canonical Bundle}

In order to find the Picard--Fuchs operator, a global section of the
canonical bundle is needed. In the case of a complete intersection,
one could simply have used the dual of $\bigwedge_j dp_j$ or its
residue form $\Res\bigwedge_i dx_i/\prod_j p_j$. The pfaffian variety,
however, is not a complete intersection. For $p_i$ the pfaffian of $N$
with row and column $i$ removed, the polynomials $p_{\mu_0},
p_{\mu_1}, p_{\mu_2}$, $\mu_i$ a permutation of $\Z_7$, give a
complete intersection whereever the submatrix
$N_{\mu_3\mu_4\mu_5\mu_6}$ of $N$ containing rows and columns $\mu_3$,
$\mu_4$, $\mu_5$, and $\mu_6$ has rank 4: ie., its pfaffian
$\Pf(N_{[\mu_3\mu_4\mu_5\mu_6]})$ is different from zero.  This
follows from $N\cdot P=0$.  Hence,
\begin{equation}\label{EQdualform}
(-1)^\mu \frac{dp_{\mu_0}\wedge dp_{\mu_1}\wedge
dp_{\mu_2}}{\Pf(N_{[\mu_3\mu_4\mu_5\mu_6]})}
\end{equation}
gives a global section of $\SO_{\P^6}(7) \otimes_{\SO_{\P^6}}
\bigwedge^3{\cal N}_X^{\vee}\cong{\omega_X^\circ}^\vee$. As
$\omega_X^\circ\cong\SO_X$, this section must be non-vanishing and
independent of $\mu$.\footnote{That is, independent of $\mu$ up to a
  constant which proves to be $(-1)^\mu$: checked with
  Maple.} Hence, the dual section in $\omega_X^\circ$ is
non-vanishing.  For smooth varieties, the canonical and dualizing
sheafs are identical, $\omega=\omega^\circ$, so we get:

\begin{prop}\label{prop:can}
On the varieties $X_A$, we have a global section of the
dualizing sheaf given by\footnote{For convenience, $k$-forms should contain
the coefficient $(2\pi i)^{-k}$. This places the closed forms in the
integral cohomology.}
\begin{equation}
\Omega=\frac{(-1)^\mu(2\pi i)^3\Pf(N_{[\mu_3\mu_4\mu_5\mu_6]})\Omega_0}
 {dp_{\mu_0}\wedge dp_{\mu_1}\wedge dp_{\mu_2}}
=\Res\frac{(-1)^\mu\Pf(N_{[\mu_3\mu_4\mu_5\mu_6]})\Omega_0}
 {p_{\mu_0}p_{\mu_1}p_{\mu_2}}
\end{equation}
where $\Omega_0$ is the global section of
$\omega_{\P^6}(7)\cong\SO_{\P^6}$ given by
\begin{equation}
\Omega_0=\frac{x_0^7}{(2\pi i)^6}\cdot
 \bigwedge_{i=1}^6 d\left(\frac{x_i}{x_0}\right).
\end{equation}
The general $X_A$ is smooth, making
$\Omega$ a global section of the canonical bundle.
\end{prop}

Actually, $\Pf(N_{[\mu_3\mu_4\mu_5\mu_6]})\not=0$ specifies the
appropriate component of $p_{\mu_0}=p_{\mu_1}=p_{\mu_2}=0$.

\section{The Orbifold Construction}

There are maps $\sigma:e_i\mapsto e_{i+1}$ and $\tau:e_i\mapsto
e_iw^i$ where $w=e^{2\pi i/7}$ and $(e_i)$ is a fixed basis for $E$,
forming a group action on $E$. The commutator is multiplication with a
constant, so in the projective setting, these two maps commute giving
an abelian $7\times7$-group $G$: eg., it gives an action on $\P(E\wedge
E)$.  We take the family of 6-planes in $\P(E\wedge E)\cong\P^{20}$
such that these maps restrict to them: ie.,
$\mbox{Span}\{\sum_{i+j\equiv k}y_{i-j} e_i\wedge e_j\}_{k\in\Z_7}$ or
in matrix representation, $N=[x_{i+j}y_{i-j}]_{i,j\in\Z_7}$, where we
take $x_i$ to be coordinates on $\P^6$ and $y_i+y_{-i}=0$.  This gives
a $\P^2$-family of 6-planes as parametrized by $[y_1:y_2:y_3]$, thus
defining a $\P^2$-subfamily of $X_A$. These have double-points at the
49 points $[x_i]_{i\in\Z_7}\in\{g([y_i]_{i\in\Z_7})|g\in G\}$.

For any 7-subgroup of $G$, there are
7 fixed-points in $\P^6$ under its action, and three lines in the
$\P^2$ parameter space such that these fixed points lie in the
corresponding varieties. We are free to choose any such subgroup, and
any of the three lines, without loss of generality, as the normalizer
of $G$ acts transitively on the eight triplets of lines.

Let $H$ be the subgroup generated by $\tau$, and choose the line
$y_3=0$. We may then use the coordinate $y=y_2/y_1$ to parametrize our
$\P^1$-family. We then have a matrix $N_y$ whose rank 4 locus defines
a degree 14 dimension 3 variety $X_y\subset\P^6$. In addition to the
49 double-points, the 7 fixed-points under $\tau$ are also
double-points.  In general, these are the only singular
points.\footnote{This has been checked using Macaulay \cite{Macaulay}
for the case $y=1$.}

For $y=0$ and $y=\infty$, the variety $X_y$ decomposes into 14
distinct 3-planes intersecting on the coordinate planes.

In addition to the line-triplet we have chosen, there are seven other
equivalent line-triplets. These intersect our chosen line in 21
points: $y^{21}-289y^{14}-58y^7+1=0$. For these values of $y$, the
variety gains 7 further double-points.

Using a construction similar to that of Candelas et.\ al.
\cite{Cand}, let $M_y=\widetilde{X_y/H}$ by a minimal (canonical)
desingularization of the quotient \cite{RoanRes}.

The map $x_i\mapsto x_iw^{5i^2}$ in the normalizer has the same effect
as $y\mapsto yw$. Hence, the natural parameter is $\phi=y^7$, and the
manifold is denoted $M_\phi$.

To give a brief review of the definition (in matrix notation):

\begin{defin}
Let $N$ be the skew-symmetric matrix $[x_{i+j}y_{i-j}]_{i,j\in\Z_7}$
where $y_i+y_{-i}=0$, and $P=[p_i]$ the pfaffians of the $6\times6$
diagonal minors; denote by $X_Y$, $Y=[y_i]$, the zero locus of $P$.

For $y_3=0$, let $y=y_2/y_1$ and denote the variety $X_y$. Let
$H=\langle\tau\rangle$ be the group acting on $X_y$ by
$\tau:x_i\mapsto wx_i$. We take a minimal desingularization of
$X_y/H$, parametrize this family by $\phi=y^7$, and denote the
resulting family of threefolds $M_\phi$.
\end{defin}

Gaining and resolving double-points corresponds to collapsing an $S^3$
to a point and then blowing it up to a $\P^1$. This increases the
Euler-characteristic by 2, either by increasing $h^{1,1}$ and
$h^{2,2}$ by one each or by reducing $h^{1,2}$ and $h^{2,1}$ by one
each. The blow-ups are along codimension 1 surfaces going through the
double-points, and each such blow-up provides us with an extra
$(1,1)$-form. Neither of these processes, the collapsing and the
blowing up, affect the dualizing sheaf as both processes are local:
contained in a set with no codimension 1 subvariety.

The creation and resolving of the $49+7$ double-points thus increases
the Euler-characteristic to 14. The action of $H$ has 14 fixed points:
two on each $\P^1$ from the blowing up of the initial fixed-points.
These quotient singularities can be desingularized without affecting
the dualizing sheaf \cite{RoanRes}. The Euler characteristic of the
desingularized quotient is given by Roan in \cite{Roan} to be 98 using
\begin{equation}
\chi(\widetilde{V/H})=\sum_{g,h\in H}\frac{\chi(V^g\cap V^h)}{|H|}
\end{equation}
for any smooth $V$, $V^g$ the fixed-point set in $V$ of $g$, $H$
an abelian group.

To see the effect on the betti numbers, we have to go in some greater
detail as to the blow-ups. We may blow the variety up along the
surfaces
$S_i=\{x_i=x_{i-3}=x_{i+3}=x_{i-2}x_{i+2}-y^2x_{i-1}x_{i+1}=0\}$,
$i\in\Z_7$, in any order. From this, $h^{1,1}$ is increased by 7.
Hence, for $\widetilde{X_y}$ we get $h^{1,1}=h^{2,2}=8$, causing
$h^{1,2}=h^{2,1}=1$. This corresponds with the dimension of the
parameter space which must therefore be the entire moduli-space.

The surfaces $S_i$ are invariant under $\tau$, hence, dividing
$\widetilde{X_y}$ with $H$ does not change $h^{1,1}$. By \cite{Bat},
the effect of dividing with the group $H$ and resolving the
fixed-point singularities, increases $h^{1,1}$ and $h^{2,2}$ by 3 for
each fixed-point, thus making $h^{1,1}=h^{2,2}=50$.

The variety now being smooth, the trivial dualizing sheaf is again
identical to the canonical sheaf, which must therefore be trivial too.

It should be pointed out that the resolution may not be unique.
However, different resolutions will merely correspond to different
parts of the A-model (K\"ahler structure related) moduli space: eg., a
flop corresponds to changing the sign of one component of $H^{1,1}$,
thereby moving the K\"ahler cone \cite{MorLec}.

\begin{prop}
  For general $y\in\P^1$, the manifolds $M_y=\widetilde{X_y/H}$ are
  Calabi--Yau manifolds having $\chi(M_y)=98$, $h^{1,1}=h^{2,2}=50$,
  and $h^{1,2}=h^{2,1}=1$; the global section of the canonical sheaf
  inherited from $X_A$ as given by \ref{prop:can}. At the points
  $\phi=0$ and $\phi=\infty$, the variety decomposes into 14 3-planes,
  and for $1-57\phi-289\phi^2+\phi^3=0$, where $y$ lies on an
  intersection between two special lines in the $\P^2$ parameter
  space, there is an extra double-point.
\end{prop}

The families $M_y$ and $X_A$ thus look like good mirror candidates.

\section{Mirror Symmetry}

We now have a 1-parameter family of Calabi-Yau manifolds $M_\phi$ with a
global section $\Omega(\phi)$ of the canonical bundle given. By the
mirror symmetry conjecture, there is a special point in our moduli
space corresponding to the `large radius limit'. Around this point, $H^3$
should have maximally unipotent monodromy. As $M_\phi$ degenerates
into 14 3-planes for $\phi=0$ (and for $\phi=\infty$) we will start
off with this as the assumed special point.

Following Morrison \cite{Mor}, there should be Gauss--Manin flat
families of 3-cycles $\gamma_0, \gamma_1$, i.e.{} sections of
$R_3\pi_*\C$, defined in a punctured neighborhood of $\phi=0$ with
$f_i(\phi)=\int_{\gamma_i(\phi)} \Omega(\phi)$, such that $f_0$
extends across $\phi=0$ and $f_1/f_0=g+\log\phi$ where $g$ extends
across $\phi=0$. The natural coordinate $t=t(\phi)$ is then given by
$t=f_1/f_0$: ie., the complexified K\"ahler structure on the mirror is
$\omega=t\omega_0$ where $\omega_0$ is a fixed K\"ahler form (the dual
of a line). As this enters only as $\exp{\int_{\eta}\omega}$, we may
use the coordinate $q=e^t=\phi e^g$.

The curve count on the mirror is arrived at using the mirror symmetry
assumtion: that the B-model Yukawa-coupling derived from the variation
of complex structure (Hodge-structure) should be equal to the A-model
Yukawa-coupling on the mirror. The A-model Yukawa-coupling is
expressed in terms of the corresponding K\"ahler structure given by
the natural coordinate and the number of rational curves in any curve
class (ie., of any given degree) by
\begin{equation}
\kappa_{ttt}=n_0+\sum_{d=1}^\infty n_d {\textstyle\frac{d^3q^d}{1-q^d}}.
\end{equation}
The B-model Yukawa-coupling $\kappa_{ttt}$ may be defined as in
\cite{Mor}, by
\begin{equation}
\kappa_{ttt}
=\kappa_{\frac{d}{dt}\frac{d}{dt}\frac{d}{dt}}
=\kappa\cdot
 {\textstyle\frac{d}{dt}\otimes\frac{d}{dt}\otimes\frac{d}{dt}}
=\int_{M_\phi}
\hat\Omega\wedge\nabla_{\frac{d}{dt}}^3\hat\Omega
\end{equation}
with $t$ the parameter on the moduli-space, $\frac{d}{dt}$ seen as
a tangent vector on the moduli space, $\nabla$ the Gauss--Manin
connection, and $\hat\Omega=\Omega/f_0$ the normalized
canonical form. (In the following, I will write
$\nabla_u=\nabla_{\frac{d}{du}}$ for any parameter $u$.)

All of this can be determined from knowing the Picard--Fuchs equation
\cite{Mor}.
The Picard--Fuchs equation is a differential equation on the parameter
space whose solutions are $\int_{\gamma(\phi)}\Omega(\phi)$ for
$\gamma$ Gauss--Manin flat sections on $R_3\pi_*\C$: ie.,
$\gamma(\phi)=\sum_i u_i \nu_i(\phi)$ where $\nu_i(\phi)\in
H_3(M_\phi,\Z)$, $u_i\in\C$. This equation has order
4: ie., for $f=\int_\gamma\Omega$, where
$\gamma=\gamma(\phi)\in\Gamma(R_3\pi_*\C)$ is any $\nabla$-flat
section of 3-cycles, we have
\begin{equation}
\int_{\gamma(\phi)}
\sum_{i=0}^4 A_i(\phi)
(\nabla_{\phi\frac{d}{d\phi}})^i\Omega
=\sum_{i=0}^4 A_i(\phi)
D_\phi^i f(\phi)=0
\end{equation}
for $D_\phi=\phi\frac{d}{d\phi}=d/d\log\phi$ the logarithmic derivative.
Maximally unipotent monodromy around $\phi=0$ is equivalent to having
$A_i(0)=0$ for $i<4$ and $A_4(0)\not=0$.

First, I will find $\gamma_0$ and calculate $f_0$. From this, I will
determine the Picard--Fuchs equation. Knowing the Picard--Fuchs
equation, $f_1$ can be found as another special
solution. Furthermore, the Yukawa coupling, $\kappa$,
satisfies a differential equation expressed in terms of the
$A$-coefficients.

\section{The Pfaffian Quotient Near $\phi=0$}

For simplicity, all calculations are pulled back from the manifold
$M_\phi$ to the variety $X_y\subset\P^6$.  At $y=0$, the variety
$X_y\subset\P^6$ degenerates into 14 3-planes intersecting along
coordinate axes, the group $H$ acting on each 3-plane. One of these
planes is given by $x_4=x_5=x_6=0$. Let $\gamma_0(0)$ be the cycle
given on this 3-plane (minus the axes) by $|x_i/x_0|=\epsilon$ for
$i=1,2,3$. We may extend this definition by continuity to a
neighborhood of $y=0$.

Rather than working with $\gamma_0$, it is more convenient to work
with the 6-cycle $\Gamma$ on $\P^6\setminus X_y$ given by
$|x_i/x_0|=\epsilon$ for $i=1,2,3$ and $|x_i/x_0|=\delta$ for
$i=4,5,6$, and view $\Omega$ as the residue of
\begin{equation}
\Psi=\Omega \wedge \frac{dp_{\nu_0}}{2\pi ip_{\nu_0}} \wedge
\frac{dp_{\nu_1}}{2\pi ip_{\nu_1}} \wedge
\frac{dp_{\nu_2}}{2\pi ip_{\nu_2}}
=\frac{(-1)^\nu\Pf(N_{\nu_3\nu_4\nu_5\nu_6})}
 {(2\pi i)^6 p_{\nu_0}p_{\nu_1}p_{\nu_2}}
\cdot
x_0^7 \bigwedge_{i=1}^6 d\left(\frac{x_i}{x_0}\right)
\end{equation}
for any permutation $\nu$ of $0,\ldots,6$. We now get
\begin{equation}
f_0(\phi)=\int_{\gamma_0(\phi)} \Omega(\phi)=\int_\Gamma \Psi(\phi)
\end{equation}
where the last integral is over a cycle which is independent of
$\phi$.

In order to make the numerator as simple as possible, choose
$\nu_1,\nu_2,\nu_3=0,3,4$. This makes
$\Pf(N_{\nu_1\nu_2\nu_5\nu_6})=x_3x_4$. Setting $x_0=1$ for simplicity
(or writing $x_i$ for $x_i/x_0$), the integral becomes
\begin{equation}
\int_\Gamma \frac{x_3x_4}{p_0p_3p_4}
\cdot \bigwedge_{i=1}^6 \frac{dx_i}{2\pi i}
=\int_\Gamma
\frac{1}{\prod_{i=0,3,4}(1-\sum_{j=1}^4 v_{i,j})}
\cdot \bigwedge_{i=1}^6 \frac{dx_i}{2\pi i x_i}
\end{equation}
where
\begin{equation}
[v_{i,j}]_{\sstack{i=0,3,4}{j=1,\ldots,4}}=\left[ \begin{array}{cccc}
\frac{x_2x_5}{x_3x_4}\cdot y & \frac{x_4x_6}{x_3}\cdot y^2 &
\frac{x_1x_3}{x_4}\cdot y^2 & -\frac{x_1x_6}{x_3x_4}\cdot y^3\\
\noalign{\medskip}
\frac{x_1x_4}{x_2x_3}\cdot y & \frac{x_2}{x_3x_6}\cdot y^2 &
\frac{x_3x_5}{x_2x_6}\cdot y^2 & -\frac{x_5}{x_2x_3}\cdot y^3\\
\noalign{\medskip}
\frac{x_3x_6}{x_4x_5}\cdot y & \frac{x_5}{x_1x_4}\cdot y^2 &
\frac{x_2x_4}{x_1x_5}\cdot y^2 & -\frac{x_2}{x_4x_5}\cdot y^3
\end{array}\right].
\end{equation}

Taking the power expansion of the right hand fraction in terms of
$v_{i,j}$, the only terms that give a contribution are products
$v^n=\prod_{i,j}v_{i,j}^{n_{i,j}}$ that are independent of the $x_i$.
The ring of
products of $v_{i,j}$ which do not contain $x_i$ is $\C[r_i]$ where
(see appendix for description of method for finding the $r_k$)
\begin{equation}
\begin{array}{l}
r_1=v_{1,4}v_{2,3}v_{3,3}=-y^7=-\phi\\
r_2=v_{1,2}v_{2,3}v_{3,4}=-y^7=-\phi\\
r_3=v_{1,3}v_{2,4}v_{3,3}=-y^7=-\phi\\
r_4=v_{1,2}v_{2,2}v_{2,3}v_{3,1}=y^7=\phi\\
r_5=v_{1,3}v_{2,1}v_{3,2}v_{3,3}=y^7=\phi\\
r_6=v_{1,1}v_{2,1}v_{2,3}v_{3,1}v_{3,3}=y^7=\phi.
\end{array}
\end{equation}

Instead of evaluating the sum over $v^n$, we may now evaluate the sum
over $r^m$ including as weights the number of times the term $r^m=v^n$
occures. This makes the integral, using the appropriate correspondence
between $m$ and $n$,
\begin{equation}\begin{array}{rcl}
{\ds\int_\Gamma\Psi(\phi)}
&=&{\ds\sum_{\smallbox{\sstack{(m_i)\in\N_0^6}{m=\sum_i m_i}}}}
 (-1)^{m_1+m_2+m_3} \phi^{m}
\prod_i {n_i\choose n_{i,1},n_{i,2},n_{i,3},n_{i,4}}\\
\noalign{\medskip}
&=&{\ds\sum_{\smallbox{\sstack{m_1,m_6,u_1,u_2\in\N_0}
 {m=m_1+m_6+u_1+u_2}}}}\;
(-1)^{m_1} \phi^{m}\cdot
\frac{m!}{m_1!m_6!u_1!u_2!(m-u_1)!(m-u_2)!}\\\noalign{\medskip}
&&\quad{\ds\cdot\sum_{\smallbox{m_2+m_4=u_1}}}\; (-1)^{m_2}
\frac{(m+m_4+m_6)!}{m_2!m_4!(m_4+m_6)!}
{\ds\cdot\sum_{\smallbox{m_3+m_5=u_2}}}\; (-1)^{m_3}
\frac{(m+m_5m_6)!}{m_3!m_5!(m_5+m_6)!}\\\noalign{\medskip}
&=&{\ds\sum_{\smallbox{\sstack{m_1,m_6,u_1,u_2\in\N_0}
 {m=m_1+m_6+u_1+u_2}}}}\;
(-1)^{m_1}\phi^{m}\cdot
{m\choose u_1}^2 {m\choose u_2}^2 {m+m_6\choose m}
{m+m_6\choose m_1,u_1+m_6,u_2+m_6}\\\noalign{\bigskip}
&=&\ds 1+5\phi+109\phi^2+3317\phi^3+121501\phi^4+\cdots
\end{array}\end{equation}
 
This function, $f_0$, should be a solution to a Picard--Fuchs equation
given by $\sum_{i=0}^4 A_i D_\phi f_0(\phi)=0$, where
$D_\phi=\phi\frac{d}{d\phi}$ and $A_i$ are polynomials in
$\phi$ with $A_i(0)=0$ for $i<4$. Entering general
polynomials for $A_i$, we find a solution for $\deg A_i=5$:
\begin{equation}\label{PFpfaff0}
\sum_{i=0}^4
 A_i D_\phi^i
\begin{array}[t]{cl}
=&(1-57\phi-289\phi^2+\phi^3)(\phi-3)^2 D_\phi^4\\
&+4\phi(\phi-3)(85+867\phi-149\phi^2+\phi^3) D_\phi^3\\
&+2\phi(-408-7597\phi+2353\phi^2-239\phi^3+3\phi^4) D_\phi^2\\
&+2\phi(-153-4773\phi+675\phi^2-87\phi^3+2\phi^4) D_\phi\\
&+\phi(-45-2166\phi+12\phi^2-26\phi^3+\phi^4).
\end{array}\end{equation}
This is the so called Picard--Fuchs operator.

Solving for $f_1(\phi)=f_0(\phi)\cdot(g(\phi)+\log\phi)$, we get
$g(\phi)=\alpha+14\phi+287\phi^2+\cdots$, where $\alpha$ is a
constant. The natural coordinate is $t=g(\phi)+\log\phi$ or
$q=e^t=c_2(\phi+14\phi^2+385\phi^3+\cdots)$ where
$c_2=e^\alpha$.

We then calculate the Yukawa coupling. This is a symmetric 3-tensor on
the parameter space, $\P^1$, which will be globaly defined but with
poles. The Yukawa coupling is given by (\cite{Mor},\cite{Cand})
\begin{equation}\begin{array}{rcl}
\ds\kappa_{ttt}
&=&\ds\left(\frac{d\log\phi}{dt}\right)^3 \kappa_{\log\phi\log\phi\log\phi}\\
&=&\ds\left(\frac{d\log\phi}{dt}\right)^3
 \int_{M_\phi} \hat\Omega\wedge
  \nabla_{\phi\frac{d}{d\phi}}^3\hat\Omega\\
&=&\ds\left(\frac{d\log\phi}{dt}\right)^3
 \frac{1}{f_0(\phi)^2}
 \int_{M_\phi} \Omega\wedge
  \nabla_{\phi\frac{d}{d\phi}}^3\Omega.
\end{array}\end{equation}
To move $f_0$ to outside the differential, we use Griffiths
transversality property which implies that
$\Omega\wedge\nabla^i\Omega=0$ for $i<3$.
 
The term
$\int_{M_\phi}\Omega\wedge\nabla_{\phi\frac{d}{d\phi}}^3\Omega$
satisfies a differential equation \cite{Mor}:
\begin{equation}
\phi\frac{d}{d\phi}\log\left(
 \int_{M_\phi}\Omega\wedge\nabla_{\phi\frac{d}{d\phi}}^3\Omega
\right)
=-\frac{A_3}{2A_4}.
\end{equation}
This gives us
\begin{equation}\label{eq:k}
\int_{M_\phi}\Omega\wedge\nabla_{\phi\frac{d}{d\phi}}^3\Omega
=\frac{c_1(3-\phi)}{1-57\phi-289\phi^2+\phi^3}
\end{equation}
for some constant $c_1$.  The denominator may be seen to have zeros at
three points in the parameter space.  These are the points where the
manifold has singularities: where our particular special line in the
bigger parameter space $\P^2$ intersects other special lines, and,
hence, has an additional double point comming from the seven extra
double points on $X_y$.

The final step is to express $\kappa_{ttt}$ in terms of $q$. Using the
power series expansion $q=q(\phi)$ and its inverse series giving
$\phi=\phi(q)$, and $\frac{d\log\phi}{dt}=\frac{q}{\phi}\frac{d\phi}{dq}$, we
may express $\kappa_{ttt}$ as
\begin{equation}\begin{array}{rcl}
\ds\kappa_{ttt}
&=&\ds\left(\frac{q}{\phi(q)}\frac{d}{dq}\phi(q)\right)^3
 \frac{1}{f_0(\phi(q))^2}\cdot
 \frac{c_1(3-\phi)}{1-57\phi-289\phi^2+\phi^3}\\\noalign{\medskip}
&=&c_1\left(3+14\frac{q}{c_2}+714(\frac{q}{c_2})^2
 +24584(\frac{q}{c_2})^3+906122(\frac{q}{c_2})^4
 +\cdots\right)\\\noalign{\medskip}
&=&c_1\left(3+14\qcterm{1}+\frac{714-14c_2}{8}\qcterm{2}
 +\cdots\right)
\end{array}\end{equation}
In order that there be only non-negative integer coefficients in the
last line, we set $c_2=1$.  Putting $c_1=2m$, we get
\begin{equation}
\textstyle
\kappa_{ttt}=m\cdot\left(6+28\qterm{1}+175\qterm{2}+1820\qterm{3}
 +28294\qterm{4}+\cdots\right).
\end{equation}
The actual value of $m$ cannot be seen from this series alone.
However, $m$ is supposed to have a fixed value as determined by the
value of the Yukawa coupling.

\begin{prop}
The manifold $M_\phi$ has maximally unipotent monodromy around
$\phi=0$, the Picard--Fuchs equation is given by (\ref{PFpfaff0}).
Assuming $c_2=1$ and $c_1=2m$, the mirror has degree $6m$, and the
rational curve count is $28m$ lines, $175m$ conics, $1820m$ cubics,
etc.
\end{prop}

As the general $X_A$ that was initially assumed to be the mirror, has
degree 14, and the first term of the $q$-series of $\kappa_{ttt}$
gives the degree of the mirror to be a multiplum
of 6, this cannot be the case. However, the point $\phi=\infty$
remains to be checked.

There is another striking observation:\footnote{This observation was
  made by Duco van Straten.} the Picard--Fuchs equation is exactly the
same as for the A-model of $G(2,7)\subset\P^{20}$ intersected by a
general $\P^{13}$ \cite{Straten}. In this case, $m=7$.

\section{The Pfaffian Quotient Near $\phi=\infty$}

Initially, the Picard--Fuchs equation seems to be regular at infinity,
which would be most surprising as $M_\infty$ degenerates into
14 3-planes just like $M_0$. However, global sections of the canonical
bundle $\Gamma(\omega_{M_\phi})$ may be viewed as a line-bundle on the
parameter space, and as such it is isomorphic to $\SO_{\P^1}(1)$. To
see this, recall that the global section $\Omega$ was of degree $-7$
in $y$, hence, degree $-1$ in $\phi$. In order to get a global section
of the canonical bundle near $\phi=\infty$, one should use
$\tilde\Omega=\phi\cdot\Omega$. This modification and changing
coordinate to $\tilde\phi=1/\phi$ amounts to the change $D_\phi\mapsto
-D_{\tilde\phi}-1$ in the Picard--Fuchs operator, making it 
\begin{equation}\label{PFpfaffinf}
\sum_{i=0}^4
 \tilde A_i D_{\tilde\phi}^i
\begin{array}[t]{cl}
=&(1-289\tilde\phi-57\tilde\phi^2+\tilde\phi^3)(1-3\tilde\phi)^2
 D_{\tilde\phi}^4\\
&+4\tilde\phi(3\tilde\phi-1)
 (143+57\tilde\phi-87\tilde\phi^2+3\tilde\phi^3)
 D_{\tilde\phi}^3\\
&+2\tilde\phi(-212-473\tilde\phi+725\tilde\phi^2
              -435\tilde\phi^3+27\tilde\phi^4)
 D_{\tilde\phi}^2\\
&+2\tilde\phi(-69-481\tilde\phi+159\tilde\phi^2
              -171\tilde\phi^3+18\tilde\phi^4)
 D_{\tilde\phi}\\
&+\tilde\phi(-17-202\tilde\phi-8\tilde\phi^2
             -54\tilde\phi^3+9\tilde\phi^4).
\end{array}\end{equation}
We now see that the monodromy is maximally unipotent around $\phi=\infty$.

We may now procede as for the previous case, but calculating $\tilde f_0$
from the Picard--Fuchs equation rather than the opposite. This gives a
Yukawa-coupling in terms of $\tilde q$:
\begin{equation}
\kappa_{\tilde t\tilde t\tilde t}
=\tilde c_1(1+42\frac{\tilde q}{\tilde c_2}+
6958(\frac{\tilde q}{\tilde c_2})^2+\cdots)
\end{equation}
where $\tilde c_1=c_1=2m$: just enter $\phi=\infty$ into the
Yukawa-coupling \ref{eq:k} after multiplying with $\phi^2$ owing to
the transition to $\tilde\Omega=\phi\cdot\Omega$. Putting $\tilde
c_2=1$, we get
\begin{equation}
\kappa_{\tilde t\tilde t\tilde t}
=m(2+84\qtterm{1}+1729\qtterm{2}+83412\qtterm{3}+5908448\qtterm{4}+\cdots).
\end{equation}

\begin{prop}
  The manifold $M_\phi$ has maximally unipotent monodromy around
  $\phi=\infty$, the Picard--Fuchs equation is given by
  (\ref{PFpfaffinf}).  Assuming $\tilde c_2=1$ and $m=7$ to give the
  mirror degree $14$, the rational curve count is $588$ lines,
  $12103$ conice, $583884$ cubics, etc.
\end{prop}

The lines on the general pfaffian have been counted by Ellingsrud and
Str\o{}mme\footnote{Private communication. May be published by
Ellingsrud, Str\o{}mme and Peskine soon.}, and is 588.

\section{The Grassmannian G(2,7) Quotient}

Due to the equality between the B-model Picard--Fuchs operator at
$\phi=0$ for the pfaffian quotient and the A-model Picard--Fuchs
operator for an intersection of $G(2,7)\subset\P^{20}$ with
a general $\P^{13}$, it is natural to take a closer look at $G(2,7)$.
In particular, it is possible to perform an orbifold constuction on
this which is `dual' to that on the pfaffian.

The pfaffian quotient was constructed from an intersection between the
general pfaffian in $\P^{20}$ and a special family of 6-planes:
$\P^{6}_y$. We may take the family $\P^{13}_y$ of 13-planes in $P(E^\vee\wedge
E^\vee)\cong\P^{20}$ dual to $\P^6_y$, and take
$Y_y=G(2,7)\cap\P^{13}_y\subset\P^{20}$. Again, we have a group action
by $\tau:x_{i,j}\mapsto x_{i,j}w^{i+j}$ which restricts to this
intersection, and the natural coordinate being $\phi=y^7$. The
$\tau$-fixed points, $e_i\wedge e_{i+3}$, are double-point
singularities, as are the images under $\tau$ of
$(e_{i+1}-e_{i-1})\wedge((e_{i+3}-e_{i-3})+y\cdot(e_{i-2}-e_{i+2}))$. Let
$W_y$ be the desingularized quotient $\widetilde{Y_y/\tau}$. This is a
Calabi--Yau manifold \cite{Straten}. I will procede without going into
the desingularization as this has no impact on the B-model.

To summarize the definition (in matrix notation):

\begin{defin}
Let $U_y=[x_{i,j}]_{i,j\in\Z_7}$, $x_{i,j}+x{j,i}=0$, be the
skew-symmetric matrix with $x_{i+4,i-4}=-yx_{i+1,i-1}$. (This amounts
to specializing to the 13-planes $\P^{13}_y$ dual to the $\P^6_y$ used
for the pfaffians, and giving a specific coordinate system.) Let $Y_y$
denote the rank~2 locus of $U_y$ in $\P^{13}$. Divide this out with
the group action generated by $\tau:x_{i,j}\mapsto x_{i,j}w^{i+j}$,
take a minimal desingularization of this, parametrize the resulting
family of threefolds by $\phi=y^7$ denoting it $W_\phi$.
\end{defin}

In order to get an expression for the canonical form, we may look
at an affine piece of $G(2,7)$ given by $u_1\wedge u_2$ where
$u_i=[u_{i,j}]$, $i=1,2$, $j=0,\ldots,6$, and where
$u_{1,0}=u_{2,2}=1$, $u_{1,2}=u_{2,0}=0$. The defining equations then
become
\begin{equation}
u_{1,i}u_{2,i+1}-u_{1,i+1}u_{2,i}=
y\cdot(u_{1,i-2}u_{2,i+3}-u_{1,i+3}u_{2,i-2}),\quad i\in\Z_7.
\end{equation}
Now, as we have a complete intersection, we may define the canonical
form $\Omega$ as the residue of
\begin{equation}
\Psi=
\frac{\bigwedge_{i=1,3,4,5,6} du_{1,i}\wedge du_{2,i}}
 {(2\pi i)^{10} \prod_{i\in\Z_7} (u_{1,i}u_{2,i+1}-u_{1,i+1}u_{2,i}-
  y\cdot(u_{1,i-2}u_{2,i+3}-u_{1,i+3}u_{2,i-2}))}.
\end{equation}

For $y=0$, the variety decomposes. One of the components may be given
in affine coordinates by $u_1\wedge u_2$ where $u_1=[1,0,0,0,0,0,0]$,
$u_2=[0,0,1,u_{2,3},u_{2,4},u_{2,5},0]$.  We may define the 3-cycle
$\gamma_0(0)$ by $|u_{2,j}|=\epsilon$ for $j=3,4,5$, and extend this to a
neighborhood: say, $|y|<\delta$.  As for the pfaffian, we
will rather use the 10-cycle $\Gamma$ in $\P^{20}$ defined by
$|u_{i,j}|=\epsilon_{i,j}$, where again $u_i=[u_{i,j}]$ with
$u_{1,0}=u_{2,2}=1$, $u_{1,2}=u_{2,0}=0$. The actual choices of $\delta$ the
$\epsilon_{i,j}$ will be made so as to make the quotients $v_{i,j}$
defined below sufficiently small, but will otherwise be of no importance.
 
We may now rewrite the residual form so as to suite our purpose of
evaluating it as a power series in $y$:
\begin{equation}
\Psi=\frac{1}{\prod_i(1-\sum_j v_{i,j})}\cdot
 \bigwedge_{i=1,3,4,5,6}\frac{du_{1,i}\wedge du_{2,i}}{(2\pi
 i)^2 u_{1,i}u_{2,i}}
\end{equation}
where
\begin{equation}\begin{array}{c}
v_{1,1}=-y\cdot\frac{u_{1,5}u_{2,3}}{u_{2,1}},\quad
v_{1,2}=y\cdot\frac{u_{1,3}u_{2,5}}{u_{2,1}}\\\noalign{\medskip}
v_{2,1}=-y\cdot\frac{u_{1,6}u_{2,4}}{u_{1,1}},\quad
v_{2,2}=y\cdot\frac{u_{1,4}u_{2,6}}{u_{1,1}}\\\noalign{\medskip}
v_{3,1}=y\cdot\frac{u_{2,5}}{u_{1,3}}\\\noalign{\medskip}
v_{4,1}=\frac{u_{1,3}u_{2,4}}{u_{1,4}u_{2,3}},\quad
v_{4,2}=y\cdot\frac{u_{1,1}u_{2,6}}{u_{1,4}u_{2,3}},\quad
v_{4,3}=-y\cdot\frac{u_{1,6}u_{2,1}}{u_{1,4}u_{2,3}}\\\noalign{\medskip}
v_{5,1}=\frac{u_{1,4}u_{2,5}}{u_{1,5}u_{2,4}},\quad
v_{5,2}=-y\cdot\frac{1}{u_{1,5}u_{2,4}}\\\noalign{\medskip}
v_{6,1}=\frac{u_{1,5}u_{2,6}}{u_{1,6}u_{2,5}},\quad
v_{6,2}=y\cdot\frac{u_{1,3}u_{2,1}}{u_{1,6}u_{2,5}},\quad
v_{6,3}=-y\cdot\frac{u_{1,1}u_{2,3}}{u_{1,6}u_{2,5}}\\\noalign{\medskip}
v_{7,1}=y\cdot\frac{u_{1,4}}{u_{2,6}}.
\end{array}\end{equation}
In order that the power series expansion converge, we need
$\sum_j|v_{i,j}|<1$. In order to obtain this, set
$\epsilon_{1,i}/\epsilon_{2,i}<\epsilon_{1,j}/\epsilon_{2,j}$ for
$3\le i<j\le6$, and $\delta$ sufficiently small.
 
If we look at the ring generated by the $v_{i,j}$, the subring of
elements that do not contain terms $u_{i,j}$ is $\C[r_i]$, where
\begin{equation}\begin{array}{c}
r_1=v_{1,1}v_{2,1}v_{3,1}v_{4,2}v_{5,2}v_{6,2}v_{7,1}
 =y^7=\phi\\\noalign{\medskip}
r_2=v_{1,2}v_{2,1}v_{3,1}v_{4,3}v_{5,2}v_{6,1}v_{6,3}v_{7,1}
 =-y^7=-\phi\\\noalign{\medskip}
r_3=v_{1,1}v_{2,1}v_{3,1}v_{4,1}v_{4,3}v_{5,1}v_{5,2}v_{6,1}v_{6,3}v_{7,1}
 =y^7=\phi\\\noalign{\medskip}
r_4=v_{1,1}v_{2,2}v_{3,1}v_{4,1}v_{4,3}v_{5,2}v_{6,3}v_{7,1}
 =-y^7=-\phi.
\end{array}\end{equation}

For any monomial $r^m=\prod_i r_i^{m_i}$, the corresponding
$u^n=\prod_{i,j}u_{i,j}^{n_{i,j}}$ appares $\prod_i{n_i\choose
n_{i,1},\ldots}$ number of times, $n_i=\sum_j n_{i,j}$. The power
series expansion for $f_0=\int_{\gamma_0}\Omega$ will then be given by
\begin{equation}\begin{array}{rcl}
{\ds\int_\Gamma\Psi}
&=&{\ds\sum_{\smallbox{\sstack{(m_i)\in\N_0^4}{m=\sum_i m_i}}}
 (-1)^{m_2+m_4}\phi^m\cdot\prod_{i=1}^7}
 {n_i\choose n_{i,1},\ldots}\\\noalign{\medskip}
&=&{\ds\sum_{\smallbox{\sstack{(m_i)\in\N_0^4}{m=\sum_i m_i}}}
 (-1)^{m_2+m_4}\phi^m\cdot}
 {m\choose m_2}{m\choose m4}{m+m_3\choose m}\\\noalign{\medskip}
&&\quad\quad
 \cdot{m+m_2+m_3\choose m_1,m_2+m_3,m_2+m_3+m_4}
 {m+m_3+m_4\choose m_1,m_3+m_4,m_2+m_3+m_4}\\\noalign{\medskip}
&=&\ds 1+5\phi+109\phi^2+3317\phi^3+121501\phi^4+4954505\phi^5+\cdots
\end{array}\end{equation}
which may be recognized as exactly the same series as for the pfaffian
quotient. Hence, the Picard--Fuchs operator etc.\ all become the same
as for the pfaffian quotient.
 
The global sections of the canonical sheaf again forms a
$\SO_{\P^1}(1)$ line-bundle on the $\P^1$ parameter space. Hence, this
grassmannian quotient has the same Picard--Fuchs operator at
$\phi=\infty$ as the pfaffian quotient.

\begin{prop}
The B-models of $M_y$ and $W_y$ have the same Picard--Fuchs
operator. Hence, the Yukawa-coupling may at most differ by a factor.
\end{prop}

Of course, it is natural to conjecture that the Yukawa-couplings are
equal, making the B-models isomorphic.

\section{Comments on the Results}

Apparantly, there is a strong relation between the varieties defined
by the pfaffians and the grassmannian $G(2,7)$. The B-models of the
$M_y$ and $W_y$ are isomorphic, and according to mirror symmetry and
assuming that we actually have the mirrors, the A-models of the
general pfaffian and general $G(2,7)$ sections should also be
isomorphic, and vice versa. It may of course be possible that we have
found models with the same B-model but different A-models, in which
case they would not be mirrors.

Assuming that we actually have mirror symmetry, it would appare that
varying the complex structure on $M_y$ or $W_y$ leads to a transition
from the K\"ahler structure on the pfaffian section $X_A$ to that of
the grassmannian section $Y_A$.

\begin{conj}
The pairs $M_y+W_y$ is the mirror family of $X_A+Y_A$ where $M_y$ and
$W_y$ (resp.\ $X_A$ and $Y_A$) form different parts of the A-model
(K\"ahler) moduli space.
\end{conj}

It is worth noting that the varieties $X_A$ and $Y_A$ cannot be
birationally equivalent. As $h^{1,1}=1$, this has a unique positive
integral generator (the dual of a line); if birational, these two must
correspond up to a rational factor. Integrating the third power of
this over the variety gives the degree; the ratio of the degrees would
then be the third power of a rational number, which is not possible
for $42/14=3$.

There may be rational 3-to-1 map though.

\appendix

\section{Finding Generators of Subring}
\label{app:subring}

Assume that we have a list of variables $x_i$, $i=1,\ldots,n$, and
Laurent-monomials $v_j=\alpha_j \prod_i x_i^{a_{j,i}}$,
$j=1,\ldots,m$, with $a_{j,i}\in\Z$. We wish to find $r_k=\prod_j
v_j^{b_{k,j}}$ such that $r_k(x)$ is independent of $x_i$ and
generates the ring of polynomials in $v_j(x)$ independent of $x_i$ (or
some extension of this ring).

An optimistic approach is simply the find a set of linearly
independent vectors with integer coefficients generating the kernel of
the matrix $A=[a_{j,i}]:\C^m\pil\C^n$. In the nicest cases, in
particular in the two cases that we are treating, one may even find
such vectors with non-negative integer coefficients. If these vectors
are $b_k=[b_{k,j}]_j$, $k=1,\ldots,m-n$, define $r_k=\prod_j
v_j^{b_k,j}$.

More generally, there is a risk that some of the $r_k$ will not be
monomials, but Laurent monomials: some $b_{k,j}$ will be negative.
These can still be used as generators, but in the sum over monomials
in $r_k$, only those which are monomials in $v_j$, ie.{} without
negative powers of $v_j$, are considered.

\section*{Acknowledgements}

I wish to thank my supervisor Geir Ellingsrud, Stein Arild Str\o{}mme,
Duco van Straten, Victor Batyrev, Sheldon Katz, Erik Tj\o{}tta, and
all others who were at the Mittag-Leffler institute 1996--97.

\raggedright

\end{document}